\documentclass[psamsfonts, a4paper, twoside]{proc-l}
\usepackage{amssymb, url}

\issueinfo{00}
{}
{}
{2001}

\copyrightinfo{2005}{American Mathematical Society}

\newtheorem{theorem}{Theorem}[section]
\newtheorem{lemma}[theorem]{Lemma}
\newtheorem{corollary}[theorem]{Corollary}

\newtheorem*{carrier theorem}{Carrier Theorem}
\newtheorem*{weak carrier theorem}{Weak Carrier Theorem}

\theoremstyle{definition}
\newtheorem*{definition}{Definition}

\theoremstyle{remark}

\newtheorem*{remark*}{Remark}

\numberwithin{equation}{section}

 \DeclareMathOperator{\st}{st}
\DeclareMathOperator{\bst}{bst} \DeclareMathOperator{\dom}{dom}
\DeclareMathOperator{\Cl}{Cl}

\newcommand{\F}[1]{{\mathcal{F}_{#1}}}
\newcommand{\G}[1]{{\mathcal{G}_{#1}}}
\newcommand{\E}[1]{{\mathcal{E}_{#1}}}

\newcommand{\SF}{{\mathcal{S}_\F{}}}
\newcommand{\BF}{{\mathcal{B}_\F{}}}

\begin{document}

\bibliographystyle{abbrv}

\title{Carrier and Nerve Theorems\\ in the Extension Theory}

\author{Andrzej Nag\'orko}
\address{Institute of Mathematics, Polish Academy of
Sciences, ul. \'Sniadeckich 8, 00-956 Warszawa, Poland.}
\email{amn@impan.gov.pl}
\thanks{The author is grateful to Professor Henryk Toru\'nczyk for
his advice while preparing this paper.}

\date{May 15, 2005 and, in revised form, August 19, 2005.}
\commby{Alexander Dranishnikov}

\begin{abstract}
  We show that a regular cover of a general topological space
  provides structure similar to a triangulation. In this general
  setting we define analogues of simplicial maps and prove their
  existence and uniqueness up to homotopy. As an application we
  give simple proofs of sharpened versions of nerve theorems of
  K.~Borsuk and A.~Weil, which state that the nerve of a regular
  cover is homotopy equivalent to the underlying space.

  Next we prove a nerve theorem for a class of spaces with
  uniformly bounded extension dimension. In particular we prove
  that the canonical map from a separable metric $n$-dimensional
  space into the nerve of its weakly regular open cover induces
  isomorphisms on homotopy groups of dimensions less than $n$.
\end{abstract}

\maketitle

\section{Introduction}

An efficient way to investigate properties of a topological space
is to divide it into pieces and examine how they are glued
together. We show how to divide a general topological space and
endow it with a structure that resembles a triangulation. We
employ this analogy to sharpen prior results of homotopy theory,
known as carrier and nerve theorems.

In order to divide a space that belongs to a class~$\mathcal{C}$
of topological spaces we must decide what a piece is. We want it
to resemble a simplex as much as possible. As a simplex is an
archetype of an absolute extensor, the choice of absolute
extensors for~$\mathcal{C}$ as pieces is quite natural.

\begin{definition}
  A space~$Y$ is an \emph{absolute extensor} for a space~$X$ if
  each map from a closed subset of~$X$ into~$Y$ extends over the
  entire space~$X$. The class of absolute extensors for all spaces
  from a class~$\mathcal{C}$ is denoted by~$AE(\mathcal{C})$. We
  write $AE(X)$ to abbreviate $AE(\{ X \})$.
\end{definition}

The following defines regular covers that endow a space with
structures similar to triangulations. Recall that a cover is
locally finite dimensional if its nerve is such.

\begin{definition}
  Let~$\mathcal{C}$ be a class of topological spaces. We say that
  a cover is a $\mathcal{C}$-cover if the intersection of each
  non-empty collection of its elements belongs to~$\mathcal{C}$. A
  locally finite $AE(\mathcal{C})$-cover that is either closed and
  locally finite dimensional or is open is said to be
  \emph{regular} for the class~$\mathcal{C}$.
\end{definition}

  Examples of regular covers include a locally finite cover of a
Euclidean space by its open balls and a cover of a finite
simplicial complex by its simplices.

  The main result of this paper states that the nerve of a regular
cover reflects the homotopy structure of the underlying space. The
result is extended to classes of spaces with bounded extension
dimension, where appropriate theory of $[L]$-homotopy equivalences
is used. Apart from generalizations, the unified approach
developed in this paper delivers simplifications of proofs of
earlier known theorems. In particular our carrier theorems imply
nerve theorems of K.~Borsuk~\cite[p. 234]{borsuk1948} and
A.~Weil~\cite[p. 141]{weil1952}, while
Theorem~\ref{thm:n-dimensional nerve theorem} generalizes an
$n$-connectivity nerve theorem recently proved by
A.~Bj\"orner~\cite{bjorner2003}.

More specialized definitions of covers that by our definition are regular were
already used in the literature. K.~Kawamura gave characterizations of infinite
dimensional manifolds in terms of partitions~\cite{kawamura1995}. G.~de~Rham
used covers by convex subsets to define simple homotopy types of Riemannian
manifolds~\cite{derham1970}. Finally nerve theorems are often used as a bridge
between combinatorics and topology~\cite{bjorner1985}.  Our generalization is
motivated by the application of regular covers in the recent proof of
characterization and rigidity theorems for N\"obeling
manifolds~\cite{nagorko2006a}. In~\cite{chigogidze2002} A.~Chigogidze
conjectured that analogous characterizations hold for universal spaces for
extension dimension.  At the end of section~\ref{sec:nerve theorems} we give a
nerve theorem for spaces of bounded extension dimension as a first step in a
program to prove these conjectures.

\section{Carrier theorems}
\label{sec:carrier theorems} How to extend a partial map from a
subcomplex to the entire CW complex? If each map from the boundary
of a Euclidean ball into the codomain extends over the ball, then
the answer is easy: order cells by inclusion and construct an
extension inductively. But the asphericity of the codomain
(vanishing of all its homotopy groups) is a rare luxury. The same
technique would work though if we were able to restrict ranges of
the map on individual cells of the CW complex to aspherical
subspaces of the codomain. This idea leads to the notion of a
carrier and to the aspherical carrier theorem~\cite[II \S
9]{lundell1969}. We generalize this notion to arbitrary spaces
using a cover to replace the cell structure in the domain.

\begin{definition}
  A \emph{carrier} is a function $C \colon \F{} \to \G{}$ from a
  cover~$\F{}$ of a space~$X$ into a collection~$\G{}$ of subsets
  of a topological space such that for each $\mathcal{A} \subset
  \F{}$ if $\bigcap \mathcal{A} \neq \emptyset$, then $\bigcap_{A
  \in \mathcal{A}} C(A) \neq \emptyset$. We say that a map~$f$ is
  \emph{carried by~$C$} if it is defined on a closed subset of~$X$ and
  $f(F) \subset C(F)$ for each $F \in \F{}$ (we write $f(F)$ for
  $f(F \cap \dom f)$).
\end{definition}

\begin{carrier theorem}
  Assume that $C \colon \F{} \to \G{}$ is a carrier such
  that~$\F{}$ is a closed cover of a space~$X$ and~$\G{}$ is an
  $AE(X)$-cover of another space. If~$\F{}$ is locally finite and
  locally finitely dimensional, then each map carried by~$C$
  extends to a map of the entire space~$X$, also carried by~$C$.
\end{carrier theorem}

Special cases of the Carrier Theorem follow from Michael's
Selection Theorem, as the multivalued map given by the formula
$F(x) = \bigcap_{F \owns x} C(F)$ is lower semi-continuous.

The proof of the Weak Carrier Theorem (later in this section) is
similar to the proof of the Carrier Theorem and can be read
simultaneously. The key difference is in the definition and an
order of sets~$\delta_\gamma$. The reader may compare both cases
for a cover~$\F{}$ that consists of three sets with non-empty
intersection.

\begin{proof}
  Let $f_0$ be a map carried by $C$ and let $A = \dom f_0$. Let
  $\{ \E{\gamma} \}_{0 < \gamma < \Gamma}$ be a transfinite
  sequence of all subcollections of~$\F{}$ with non-empty
  intersections, such that the sequence $\{ \delta_\gamma =
  \bigcap \E{\gamma}\}_{0 < \gamma < \Gamma}$ is non-decreasing in
  the order by inclusion. Its existence is guaranteed by the
  assumption of local finite dimensionality of~$\F{}$. Let
  $C_\gamma = \bigcap_{E \in \E{\gamma}} C(E)$ and observe that an
  arbitrary map~$f$ is carried by~$C$ if and only if
  $f(\delta_\gamma) \subset C_\gamma$ for each $0 < \gamma <
  \Gamma$. Let~$\delta_\Gamma = \emptyset$. We shall construct a
  transfinite sequence of maps $\{ f_\gamma \colon A \cup
  \bigcup_{0 < \alpha \leq \gamma} \delta_\alpha \to Y \}_{\gamma
  \leq \Gamma}$ such that~$f_\alpha$ extends~$f_\beta$ for all $0
  \leq \beta \leq \alpha \leq \Gamma$ and $f_\gamma(\delta_\gamma)
  \subset C_\gamma$ for each $0 < \gamma < \Gamma$. The
  map~$f_\Gamma$ will be an extension that we are looking for,
  since~$\bigcup_{\gamma < \Gamma} \delta_\gamma = X$.

  We proceed by transfinite induction. Fix $\gamma \leq \Gamma$
  and assume that for each~$\alpha < \gamma$ we already
  constructed~$f_\alpha$. The map $f'_\gamma = \bigcup_{\alpha <
  \gamma} f_\alpha$ is well defined, continuous and its domain is
  closed in~$X$ because maps~$f_\alpha$ agree on intersections of
  their domains and~$\F{}$ is closed and locally finite.
  If~$\gamma = \Gamma$, then~$\delta_\gamma = \emptyset$ and we may
  put $f_\Gamma = f'_\Gamma$. If $\gamma < \Gamma$, then by the
  order of $\delta_\gamma$ and by inductive assumptions~$f'_\gamma$
  maps $\delta_\gamma$ into~$C_\gamma$. The set~$C_\gamma$ is
  non-empty because~$C$ is a carrier and is an absolute extensor
  for~$\delta_\gamma$ because~$\G{}$ is an $AE(X)$-cover. So
  $f'_\gamma$ extends onto~$\delta_\gamma$ to a map~$f_\gamma$
  such that $f_\gamma(\delta_\gamma) \subset C_\gamma$ and our
  construction is finished.
\end{proof}

\begin{definition}
A cover is regular for a space~$X$ if it is regular for the
class~$\{ X \}$.
\end{definition}

\begin{corollary}
  \label{cor:carried homotopy} If a closed cover~$\G{}$ of a
  space~$Y$ is regular for $X \times [0, 1]$, then every two
  $\G{}$-close maps from $X$ into $Y$ are $\G{}$-homotopic.
  Moreover such homotopy exists with an additional property that
  if endpoints of its path lie in an element of~$\G{}$, then the
  entire path lies in it.
\end{corollary}
\begin{proof}
  Let $f$ and $g$ denote $\G{}$-close maps from $X$ into $Y$.  Let
  $\F{} = \{ F_G \}_{G \in \G{}}$ be the collection of subsets of $X
  \times [0, 1]$ defined by
  \begin{displaymath}
    F_G = (f^{-1}(G) \cap g^{-1}(G)) \times [0, 1].
  \end{displaymath}
  It is a cover of $X \times [0, 1]$ because $f$ and $g$ are
  $\G{}$-close. Define a carrier~$C \colon \F{} \to \G{}$  by the
  formula $C(F_G) = G$ and a map $F \colon X \times \{ 0, 1 \} \to
  Y$ by $F(x, 0) = f(x)$ and $F(x, 1) = g(x)$. By definition, $F$
  is carried by $C$ and by the Carrier Theorem it admits an
  extension to the entire space $X \times [0, 1]$, also carried by
  $C$. This extension is a $\G{}$-homotopy that satisfies our
  claim.
\end{proof}

In applications (see \cite{kawamura1995, bjorner2003}), covers with
slightly weaker regularity conditions are sometimes used, as in
the following definition.

\begin{definition}
  Let $\mathcal{C}$ be a class of topological spaces. We say that
  a cover is a \emph{weak} $\mathcal{C}$-cover if the union of
  each collection of its elements that has a non-empty
  intersection belongs to~$\mathcal{C}$. A locally finite weak
  $AE(\mathcal{C})$-cover that is either closed and locally finite
  dimensional or is open is said to be \emph{weakly regular} for
  the class $\mathcal{C}$.
\end{definition}

It follows from the Carrier Theorem that a regular cover is a
weakly regular cover, but the converse is not true.

\begin{definition}
  Let $\F{}$ be a cover of a space $X$ and let $C \colon \F{} \to
  \G{}$ be a carrier. We say that a map~$f$ is weakly carried
  by~$C$ if it is defined on a closed subset of~$X$ and for each
  $x \in \dom f$ there exists an $F \in \F{}$ such that $x \in F$
  and $f(x) \in C(F)$.
\end{definition}

If $\F{}$ is a cover of a space~$X$ and~$id_\F{}$ is the identity map
of~$\F{}$, then a map from~$X$ into~$X$ is weakly carried by $id_\F{}$ if and
only if it is~$\F{}$-close to the identity of~$X$. By this observation the
composition of two maps weakly carried by $id_\F{}$ does not have to be weakly
carried by it.

\begin{weak carrier theorem}
  Assume that~$C \colon \F{} \to \G{}$ is a carrier such that $\F{}$ is an
  open cover of a space~$X$ and $\G{}$ is a weak $AE(X)$-cover of another
  space. If $\F{}$ is locally finite, then each map weakly carried by~$C$
  extends to a map of the entire space~$X$, also weakly carried by~$C$.
\end{weak carrier theorem}

\begin{proof}
  Let $f_0$ be a map weakly carried by $C$ and let $A = \dom f_0$.  Let $\{
  \E{\gamma} \}_{0 < \gamma < \Gamma}$ be a transfinite sequence of all
  subcollections of elements of~$\F{}$ with non-empty intersections,
  non-decreasing in the order by inclusion.  For each $0 < \gamma < \Gamma$
  let $\delta_\gamma$ be the set of points in~$X$ that belong exactly to all
  elements of $\E{\gamma}$, that is, $\delta_\gamma = \bigcap \E{\gamma}
  \setminus \bigcup (\F{} \setminus \E{\gamma})$. Let $C_\gamma = \bigcup_{E
    \in \E{\gamma}} C(E)$ and observe that an arbitrary map $f$ is weakly
  carried by $C$ if and only if $f(\delta_\gamma) \subset C_\gamma$ for each
  $0 < \gamma < \Gamma$. Let $\delta_\Gamma = \emptyset$.  We shall construct
  a transfinite sequence of maps $\{ f_\gamma \colon A \cup \bigcup_{0 <
    \alpha \leq \gamma} \delta_\alpha \to Y \}_{\gamma \leq \Gamma}$ such that
  $f_\alpha$ extends $f_\beta$ for all $\beta \leq \alpha \leq \Gamma$ and
  $f_\gamma(\delta_\gamma) \subset C_\gamma$ for each $0 < \gamma < \Gamma$.
  The map $f_\Gamma$ will be an extension that we are looking for, since
  $\bigcup_{\gamma < \Gamma} \delta_\gamma = X$.
  
  We proceed by transfinite induction. Fix $\gamma \leq \Gamma$ and assume
  that for each $\alpha < \gamma$ we already constructed $f_\alpha$. The map
  $f'_\gamma = \bigcup_{\alpha < \gamma} f_\alpha$ is well defined and
  continuous because maps $f_\alpha$ agree on intersections of their domains
  and $\F{}$ is locally finite. By the definition of the order of sets
  $\delta_\gamma$ its domain is closed. If $\gamma = \Gamma$, then
  $\delta_\gamma = \emptyset$ and we may put $f_\Gamma = f'_\Gamma$. If
  $\gamma < \Gamma$, then by the order of $\delta_\gamma$ and inductive
  assumptions $f'_\gamma$ maps $\delta_\gamma$ into $C_\gamma$. The
  set~$C_\gamma$ is an absolute extensor for $\delta_\gamma$ because it is a
  union of elements of $\G{}$ that have non-empty intersection and $\G{}$ is a
  weak $AE(X)$-cover. So $f'_\gamma$ extends onto $\delta_\gamma$ to a map
  $f_\gamma$ such that $f_\gamma(\delta_\gamma) \subset C_\gamma$ and our
  construction is finished.
\end{proof}

\begin{corollary}
  \label{cor:weakly carried homotopy} If an open cover~$\G{}$ of a
  space~$Y$ is weakly regular for $X \times [0, 1]$, then every two
  $\G{}$-close maps from $X$ into $Y$ are $\st \G{}$-homotopic,
  where $\st \G{}$ denotes the star of~$\G{}$. Moreover such
  homotopy exists with an additional property that if endpoints of
  its path lie in an element of~$\G{}$, then the entire path lies
  in its star.
\end{corollary}

 We omit the proof as it is similar to the proof of
Corollary~\ref{cor:carried homotopy}.

\begin{remark*}
  In the Weak Carrier Theorem the assumption that~$\F{}$ is locally finite may
  be omitted if $\G{}$ is open and $X$ is paracompact. To prove this, it
  suffices to find a locally finite open refinement $\mathcal{H}$ of $\F{}$
  and a carrier $D \colon \mathcal{H} \to \F{}$ such that $f$ is weakly
  carried by $C \circ D$ and $H \subset D(H)$ for each $H \in \mathcal{H}$.
  Then, by the Weak Carrier Theorem, $f$ extends over $X$ to a map weakly
  carried by $C \circ D$, which is obviously weakly carried by $C$. Let
  $\mathcal{H}$ be any locally finite refinement of an open cover $\{ F \cap
  f^{-1}(C(F)) \colon F \in \F{} \}$ and let $D \colon \mathcal{H} \to \F{}$
  be any map such that $H \subset D(H) \cap f^{-1}(C(D(H)))$ for each $H \in
  \mathcal{H}$. Then $f$ is carried by $C \circ D$ because for each $x \in X$
  there exists $H \in \mathcal{H}$ such that $x \in H \in D(H) \cap
  f^{-1}(C(D(H)))$ so $f(x) \in (C \circ D)(H)$. By the definition $f$ is
  weakly carried by $C \circ D$ and the proof is finished.
  
  Similarly in Corollary~\ref{cor:weakly carried homotopy} the assumption of
  local finiteness of~$\G{}$ may be replaced by paracompactness of~$X$.
\end{remark*}

\section{Nerve theorems}
\label{sec:nerve theorems}

Nerve theorems give conditions under which the nerve of a cover is equivalent
to the underlying space. First examples of such theorems are attributed to
K.~Borsuk~\cite[p. 234]{borsuk1948} (for closed covers) and A.~Weil~\cite[p.
141]{weil1952} (for open covers), both for homotopy equivalences. Since then,
several generalizations were made. First generalizations by
W.~Holszty\'nski~\cite{holsztynski1964} and J.~N.~Haimov~\cite{haimov1979}
relaxed conditions on the cover.  Next weak homotopy equivalences were studied
in this context by M.~McCord~\cite{mccord1967} and weak $n$-homotopy
equivalences by A.~Bj\"orner~\cite{bjorner2003}.

Let $K$ denote an arbitrary simplicial complex. An open star~$\st
v$ of a vertex~$v \in K$ is the complement of the union of all
simplices of~$K$ that do not contain~$v$. A~barycentric star~$\bst
v$ of a vertex~$v \in K$ is the union of all simplices of the
first barycentric subdivision of~$K$ that contain~$v$.

\begin{lemma}
  \label{lem:regular stars} If the cover of a simplicial complex (with the
  metric topology) by the collection of open (or barycentric) stars of its
  vertices is locally finite (and locally finite dimensional), then it is
  regular for the class of metric spaces. If additionally the complex is
  locally countable, then the cover is regular for the class of normal spaces.
\end{lemma}

By \cite[Theorem 11.7, p. 109]{hu1965} and by \cite[Theorem 7.1, p.
43]{hu1965}, it suffices to show that each non-empty intersection of a
collection of elements of the cover is contractible. We omit details of the
proof, which is based on the observation that every point of the intersection
is connected by a line with the barycenter of centers of stars that are
intersected.

 We shall use the following notation.

\begin{definition}
  Let~$\SF$ denote the collection of open stars of vertices of the
  nerve~$N(\F{})$ of a point-finite cover~$\F{}$. Let~$\BF$ denote
  the collection of its barycentric stars.
\end{definition}

By Lemma~\ref{lem:regular stars}, covers $\SF$ and $\BF$ are
regular for the class of normal spaces. Their structure mirrors
the structure of~$\F{}$, in the sense of the following definition.

\begin{definition}
  We say that covers are \emph{isomorphic} if there exists a
  carrier from one of them onto the other, which is
  \emph{invertible} (i.e., is bijective and its inverse is a
  carrier).
\end{definition}

Observe that if the codomain of a carrier~$C$ is equal to the domain of a
carrier~$D$, then a composition of two maps carried respectively by~$C$
and~$D$ is carried by~$D \circ C$.

\begin{lemma}
  \label{lem:invertible carriers} Let $\F{}$ be a point-finite
  cover and let $v(F)$ denote the vertex of $N(\F{})$
corresponding to a set $F \in \F{}$. The functions $S \colon \F{}
\to \mathcal{S}_\F{}$ and $B \colon \F{} \to \mathcal{B}_\F{}$
defined by $S(F) = \st v(F)$ and $B(F) = \bst v(F)$ are invertible
carriers.
\end{lemma}
\begin{proof}
  Obviously $S$ and $B$ are bijections. We shall prove that $B$
  and $S^{-1}$ are carriers. Then, since a function $I \colon \BF
  \to \SF$ defined by the formula $I(\bst v) = \st v$ is a carrier
  and $S^{-1} \circ I \circ B = id_\F{}$, functions $B^{-1} =
  S^{-1} \circ I$ and $S = I \circ B$ are carriers too.

  To prove that $B$ is a carrier, assume that $\{ F_i \}$ is a
  collection of elements of $\F{}$ with non-empty intersection.
  Then the nerve $N(\F{})$ contains a simplex $\sigma$ spanned by
  vertices $\{ v(F_i) \}$ and each barycentric star $\bst v(F_i)$
  contains the barycenter of~$\sigma$. Therefore the intersection
  of $\{ B(F_i) \}$ is non-empty.
  
  To prove that $S^{-1}$ is a carrier, assume that $\{ \st v(F_i) \}$ is a
  collection of elements of $\SF$ with non-empty intersection. As an open star
  of a vertex $v(F_i)$ contains only points of simplices that contain
  $v(F_i)$, then the intersection of the set $\{ \st v(F_i) \}$ contains only
  points of simplices that contain all vertices $v(F_i)$. So the nerve
  $N(\F{})$ contains a simplex spanned by $\{ v(F_i) \}$ and the intersection
  of $\{ F_i \}$ is non-empty.
\end{proof}


Putting everything together we obtain nerve theorems for homotopy
equivalences. We state the theorem for closed covers, which generalizes a
nerve theorem by J.~N.~Haimov~\cite{haimov1979}. An analogous theorem for open
covers may also be proved. We do not state it here, as it turns out to be
equivalent to the nerve theorem by A.~Weil~\cite{weil1952}.

\begin{theorem}
  \label{thm:nerve theorem} Assume that a closed cover~$\F{}$ of a
  normal space~$X$ is regular for the class of metric spaces.  If~$\F{}$ is
  star-countable, then~$X$ and the nerve of~$\F{}$ are homotopy equivalent.
\end{theorem}

The main theorem of~\cite{haimov1979} states the same conclusion
under the additional assumption that~$\F{}$ is star-finite and~$X$
is paracompact.

\begin{proof}
  Let $B \colon \F{} \to \BF$ be an invertible carrier as defined in
  Lemma~\ref{lem:invertible carriers}. By the Carrier Theorem and
  Lemma~\ref{lem:regular stars} there exist $f \colon X \to N(\F{})$ and $g
  \colon N(\F{}) \to X$ carried by $B$ and $B^{-1}$ respectively. Then $g
  \circ f$ is carried by $B^{-1} \circ B$ so it is $\F{}$-close to $id_X$ and
  by Corollary~\ref{cor:carried homotopy}~$g$ is a homotopy inverse of~$f$.
  Analogously~$f$ is a homotopy inverse of~$g$ so~$X$ and~$N(\F{})$ are
  homotopy equivalent.
\end{proof}

We turn our attention to spaces with bounded dimension. First we
prove a nerve theorem for the class of at most $n$-dimensional
spaces. Next the case of extension dimension is studied. To avoid
anomalies of the dimension of general topological spaces {\bf from
now on all spaces are assumed to be separable metric}.

\newcommand{\A}{{\mathcal{A}}}

\begin{theorem}
  \label{thm:n-dimensional nerve theorem} If $\F{}$ is an open
  cover of a space~$X$, weakly regular for the class of at most
  $n$-dimensional spaces, then each canonical map $\varkappa
  \colon X \to N(\F{})$ induces isomorphisms on homotopy groups of
  dimensions less than $n$.
\end{theorem}

A canonical map into the nerve of a cover is a map induced by a
partition of unity subordinated to this cover, by interpreting its
values as barycentric coordinates of points in the nerve.
Kuratowski's $\varkappa$-map~\cite[p. 321]{engelking1992} is an
example of such a map.

It follows from the Excision Theorem that an open cover~$\F{}$ of
a locally $(n-1)$-connected space is weakly regular for the class
of at most $n$-dimensional spaces if and only if an intersection
of each collection $\A \subset \F{}$ is $(n - \left| \A
\right|)$-connected. This relates our theorem to Bj\"orner's nerve
theorem~\cite{bjorner2003}.

\begin{proof}
  For each subcollection $\A$ of $\F{}$ with non-empty
  intersection let~$v(\A)$ denote the simplex of the
  nerve~$N(\F{})$ spanned by vertices~$\{ v(A) \}_{A \in \A}$.

  Let $\lambda \colon N(\F{})^{(n)} \to X$ denote a map from the
  $n$-dimensional skeleton of $N(\F{})$ into $X$, weakly carried
  by the carrier $S^{-1}$ defined in Lemma~\ref{lem:invertible
  carriers}. Fix $k < n$. For each map $\varphi \colon S^k \to
  N(\F{})$ pick a map $\varphi' \colon S^k \to N(\F{})^{(n)}$,
  homotopic to $\varphi$, such that if $\varphi(x) \in \sigma$,
  then $\varphi'(x) \in \sigma^{(n)}$ for each simplex $\sigma$ of
  $N(\F{})$. It exists by the Cellular Approximation Theorem.

  To prove that $\varkappa$ induces an epimorphism on $k$th
  homotopy groups we shall observe that for each $\varphi \colon
  S^k \to N(\F{})$ the map $\varkappa \circ (\lambda \circ
  \varphi')$ is $\SF$-close to $\varphi$, so by
  Corollary~\ref{cor:weakly carried homotopy} they are homotopic.
  Fix $x \in S^k$. From the definition of a weakly carried map
  there exists $F \in \F{}$ such that $\varphi'(x) \in \st v(F)$
  and $\lambda(\varphi'(x)) \in F$. From the definition of
  $\varphi'$ we have $\varphi(x) \in \st v(F)$. From the
  definition of $\varkappa$ we have $\varkappa(F) \subset \st
  v(F)$. Then $\varkappa(\lambda(\varphi'(x))) \in \st v(F)$ and
  we are done.

  To prove that $\varkappa$ induces a monomorphism on $k$th
  homotopy groups we shall observe that for each $\psi \colon S^k
  \to X$ the map $\lambda \circ ((\varkappa \circ \psi)')$ is
  $\F{}$-close to $\psi$, so by Corollary~\ref{cor:weakly carried
  homotopy} they are homotopic, as dimension of $S^k \times [0,
  1]$ is at most $n$. Fix $x \in S^k$ and let $\A = \{ F \in \F{}
  \colon \psi(x) \in F \}$. Then $(\varkappa \circ \psi)(x) \in
  v(\A)$ so $(\varkappa \circ \psi)'(x) \in v(\A)^{(n)}$. By the
  definition of $\SF$, $v(\A) \cap \st v(F) \neq \emptyset$ if and
  only if $F \in \A$, so $\lambda \circ (\varkappa \circ \psi)'(x)
  \in F$ for some $F \in \A$. Therefore if $\varkappa \circ \psi$
  is null-homotopic, then $\lambda \circ ((\varkappa \circ \psi)')$
  and so is~$\psi$.
\end{proof}

\begin{remark*}
  Let~$\sigma(\A) = \Cl(\bigcap \A) \setminus \bigcup (\F{}
  \setminus \A)$ be the closure of the set of points in~$X$ that
  belong exactly to elements of $\A$. Then the map $K(\sigma(\A))
  = v(\A)$, where $\A$ runs over all subcollections of $\F{}$ that
  have non-empty intersections, is a carrier. To prove that,
  observe that if $x \in \bigcap_i \sigma(\A_i)$ then $\A = \{ F
  \in \F{} \colon x \in F \} \subset \A_i$ for each $i$, so
  $\bigcap_i v(\A_i) \supset v(\A) \neq \emptyset$. Every
  canonical map into the nerve of a cover is carried by the
  carrier $K$.
\end{remark*}

We finish with a generalization of Theorem~\ref{thm:n-dimensional
nerve theorem} to a class of spaces with uniformly bounded
extension dimension. The survey~\cite{chigogidze2002} is a good
source of information about notions of the extension dimension and
the $[L]$-homotopy; here we will only recall basic definitions.

\begin{definition}
  Let $L$ be an arbitrary CW complex. A space $X$ has extension dimension less
  than or equal to $[L]$ if $L$ is an absolute extensor for $X$. The class of
  absolute extensors for at most $[L]$-dimensional spaces is denoted
  by~$AE[L]$.

  We say that maps $f, g \colon X \to Y$ are $[L]$-homotopic if
  for each at most $[L]$-dimensional space $Z$, each pair $A$, $B$
  of disjoint closed subsets of $Z$ and each map $h \colon Z \to
  X$ there exists an extension of a map $f \circ h_{|A} \cup g
  \circ h_{|B}$ to the entire space~$Z$.
\end{definition}

We shall need the following analogue of Corollary~\ref{cor:weakly
carried homotopy}.

\begin{lemma}
  \label{lem:carried l-homotopy}
  If $\G{}$ is a locally finite closed weak $AE[L]$-cover of a
  separable metric space, then every two $\G{}$-close maps are
  $[L]$-homotopic.
\end{lemma}
\begin{proof}
  Name the maps considered by $f,g \colon X \to Y$, where $Y$
  denotes the space covered by~$\G{}$. We are going to prove that
  $f$ and $g$ are $[L]$-homotopic directly from the definition.
  Let $Z$ and $h$ be as in the definition of the $[L]$-homotopy.

  Let $\E{} = \{ E_G \}_{G \in \G{}}$ be the collection of subsets
  of $Z$ defined by
  \begin{displaymath}
    E_G = h^{-1}(f^{-1}(G) \cap g^{-1}(G)).
  \end{displaymath}
  It is a cover of $Z$ because $f$ and $g$ are $\G{}$-close. The
  map $C \colon \E{} \to \G{}$ given by the formula $C(F_G) = G$
  carries $\F{}$ into $\G{}$. By the definition $H = f \circ
  h_{|A} \cup g \circ h_{|B}$ is carried by $C$.
  
  Since $\G{}$ is a locally finite, closed cover, so is $\E{}$.  Then
  by~\cite[p. 393]{engelking1989}, since $Z$ is separable metric, there is a
  locally finite open cover $\F{} = \{ F_E \}_{E \in \E{}}$ of $Z$ such that
  $E \subset F_E$ for each $E \in \E{}$ and the function $F_E \mapsto E$ is an
  invertible carrier of $\F{}$ onto $\E{}$. Let $D \colon \F{} \to \E{}$ be a
  function such that $F \supset D(F)$ for each $F \in \F{}$. Then the identity
  map $id_Z$ is weakly carried by $D$, and $H \circ id_Z$ is weakly carried by
  $C \circ D$. By the Weak Carrier Theorem, the map $H$ extends to the entire
  space $Z$ and the proof is finished.
\end{proof}

The following is a nerve theorem for $[L]$-homotopy.

\begin{definition}
  We say that topological spaces $X$ and $Y$ are $[L]$-homotopy
  equivalent if there exist mappings $f \colon X \to Y$ and $g
  \colon Y \to X$ such that $g \circ f$ and $f \circ g$ are
  $[L]$-homotopic to $id_X$ and $id_Y$ respectively.
\end{definition}

\begin{theorem}
  \label{thm:L-dimensional nerve theorem}
  If two at most $[L]$-dimensional spaces have isomorphic closed
  covers, regular for the class of at most $[L]$-dimensional
  spaces, then they are $[L]$-homotopy equivalent.
\end{theorem}

We omit the proof, which is a straightforward application of the
Carrier Theorem and Lemma~\ref{lem:carried l-homotopy}, similar to
the proof of Theorem~\ref{thm:nerve theorem}.

\begin{remark*}
Observe that an $[S^{n-1}]$-ho\-mo\-to\-py is J.~H.~C.~Whitehead's
$n$-homotopy, and by Whitehead's characterization a map between
locally $(n-1)$-connected $n$-di\-men\-sio\-nal spaces is an
$[S^{n-1}]$-homotopy equivalence if and only if it induces
isomorphisms on homotopy groups of dimensions less than~$n$. This
connects theorem~\ref{thm:n-dimensional nerve theorem} with
Theorem~\ref{thm:L-dimensional nerve theorem}.
\end{remark*}

\end{document}